\documentclass[letterpaper, 9 pt, conference]{ieeeconf} 

\title{\textbf{In-Network Linear Regression with Arbitrarily Split Data Matrices}} 
\author{Fran\c{c}ois D. C\^ot\'e, Ioannis N. Psaromiligkos, and Warren J. Gross%
\thanks{The authors are with the Department of Electrical and Computer Engineering, McGill University, Montreal, QC H3A 0E9, Canada (e-mail: francois.cote@mail.mcgill.ca; yannis@ece.mcgill.ca; warren.gross@mcgill.ca).}}

\IEEEoverridecommandlockouts                              
\usepackage{tikz,pgfplots} 
\usepackage{amsmath,amssymb}  

\newtheorem{algorithm}{Algorithm}
\newtheorem{proposition}{Proposition}

\DeclareMathOperator{\prox}{prox}

\begin{document}

\maketitle

\thispagestyle{empty}
\pagestyle{empty}

\section{Introduction}

This paper addresses for the first time the problem of how networked agents can collaboratively fit a linear model when each agent only ever has an \emph{arbitrary} summand of the regression data. 

\subsection{Problem Statement}

Consider a network of agents (these can be computers, data centers, etc.), all concerned with a given process. Each agent has amassed some measurements of some features of the process, along with observations, or labels, of those measurements. Let $X$ be a matrix comprising all the measurements of the process across the network. Suppose that there are a total of $n$ real-valued measurements (examples) of each of $p$ features, so that $X \in \mathbf{R}^{n \times p}$. Let $y$ in $\mathbf{R}^n$ be a vector comprising the (real) labels of the examples. It is possible that no agent has the full data, but the data of all the agents covers $X$ and $y$. Fig.~1 highlights what part of $X$ an agent may have. Agents may have blocks of entire columns (features), as in Fig.~1(a), entire rows (examples), as in Fig.~1(b), or blocks of partitions, as in Fig.~1(c), but also arbitrary parts, possibly nonrectangular and overlapping, as in Fig.~1(d). The latter typifies a case not yet considered in the literature.  

In this paper, our goal is to have each agent linearly regress $y$ on $X$ within a residual $\ell_2$-norm of $\epsilon$, and to regularize the solution through a cost function $f\colon \mathbf{R}^p \to \mathbf{R}\cup\{+\infty\}$. In other words, each agent must solve the optimization problem  
\begin{equation*}
P_0\colon\quad\min_{\beta \in \mathbf{R}^p}\enskip f(\beta) \quad \mathrm{s.t.} \quad \lVert X\beta-y \rVert_2 \leq \epsilon.
\end{equation*}

To clarify the context, we lay down some assumptions:
\begin{enumerate}
\item $f$ is a closed proper convex function (not necessarily smooth);
\item each agent knows $f$ and $\epsilon$, and the number of  agents, say $m$, in the network;
\item the network is connected; but
\item an agent can only communicate with its neighbors in the network; though
\item no agent can divulge its part of $(X,y)$; and
\item[$\boldsymbol{\star}$)] each agent knows where in $X$ and $y$ its part lies and which subparts are repeated among other agents and how many times.
\end{enumerate}
Let us label the agents $1$ to $m$, and comment on assumption $\boldsymbol{\star}$. Because of this key assumption, agent $i$ can express its part of the data as a matrix $X_i$ and a vector $y_i$, so that among all the agents, $X=X_1+\cdots+X_m$ and $y=y_1+\cdots+y_m$. Of course, $X_i \in \mathbf{R}^{n \times p}$ and $y_i \in \mathbf{R}^{n}$, and by assumption 5, neither can be transferred between agents.

We can now formally describe our goal. It is to solve the problem
\begin{align*}
P_1\colon\quad&\!\min_{\beta \in \mathbf{R}^p} \enskip f(\beta)  \\
&\mathrm{\hspace{7.5mm} s.t.}\quad \lVert (X_1+\cdots+X_m)\beta-(y_1+\cdots+y_m) \rVert_2 \leq \epsilon,\\
& \text{abiding by assumptions 1--5.}
\end{align*}

\begin{figure}[t!]
	\delimitershortfall = 0pt 
	\def\a{0.0625} 
	\def\b{0.675}
	\def\c{0.8}
	\def\d{1.4125}
	\tikzset{
		every picture/.style={
		scale = 0.45,
    	baseline = -0.0875cm, 
    	rounded corners = 1.5pt,
    	semithick}}
	\hspace*{\fill}
	{%
	\begin{minipage}[b]{1.435cm} 
		\centering
		\centerline{
		$\left(
		\begin{tikzpicture} 
			\filldraw (-\d, \d) rectangle (-\c, -\d); 
			\draw (-\b, \d) rectangle (-\a, -\d);
			\draw (\a, \d) rectangle (\d, -\d);
		\end{tikzpicture} 
		\right)$}
		\medskip
		\centerline{\footnotesize(a)}
	\end{minipage}
	}
	\hfill
	{%
	\begin{minipage}[b]{1.435cm}
		\centering
		\centerline{
		$\left(
		\begin{tikzpicture} 
			\filldraw (-\d, \d) rectangle (\d, \c); 
			\draw (-\d, \b) rectangle (\d, \a); 
			\draw (-\d, -\a) rectangle (\d, -\d);
		\end{tikzpicture} 
		\right)$}
		\medskip
		\centerline{\footnotesize(b)}
	\end{minipage}
	}
	\hfill
	{%
	\begin{minipage}[b]{1.435cm}
		\centering
		\centerline{
		$\left(
		\begin{tikzpicture}  
			\filldraw (-\d, \d) rectangle (-\c, \c); 
			\draw (-\b, \d) rectangle (-\a, \c);
			\draw (\a, \d) rectangle (\d, \c);
			\draw (-\d, \b) rectangle (-\c, \a); 
			\draw (-\b, \b) rectangle (-\a, \a);
			\draw (\a, \b) rectangle (\d, \a);
			\draw (-\d, -\a) rectangle (-\c, -\d); 
			\draw (-\b, -\a) rectangle (-\a, -\d);
			\draw (\a, -\a) rectangle (\d, -\d);
		\end{tikzpicture} 
		\right)$}
 		\medskip
		\centerline{\footnotesize(c)}
	\end{minipage}
	}
	\hfill
	{%
	\begin{minipage}[b]{1.435cm}
		\centering
		\centerline{
		$\left(
		\begin{tikzpicture} 
			\filldraw (-\d, \d) |- (-\a, \a) |- (-\c, \b) -- (-\c, \d) -- cycle; 
			\draw (\c, \d) rectangle (\d, \c);
			\draw (-\b, \b) rectangle (\d, -\d);
			\draw (-\d, -\a) rectangle (-\c, -\d); 
			\draw (\a, \c) |- (\b, -\b) |- (\d, \a) |- (\b, \b) -- (\b, \d) -| (-\b, \c) -- cycle;
			\draw (\a, -\c) rectangle (\d, -\d);
		\end{tikzpicture} 
		\right)$}
		\medskip
		\centerline{\footnotesize(d)}
	\end{minipage}
	}
	\hspace*{\fill}	
\caption{Examples of data splitting by (a) features, (b) examples, (c) blocks of both, and (d) nonrectangular and overlapping blocks.}
\end{figure}
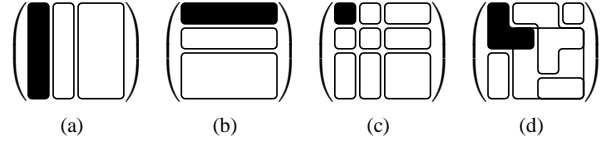

\subsection{Related Work}

To the best of our knowledge, $P_1$ has never been posed in such a general form. Many works have considered special cases of the problem, with either specific ways of splitting the data among the agents or specific network topologies. Our work is most related to a combination of the notions put forth by Mota et al.~\cite{1}--\cite{3} and Parikh and Boyd~\cite{4}, the former considering feature and example splittings of $X$ (see Fig.~1(a) and 1(b)), while the latter considering block splitting (see Fig.~1(c)) tied to a specific network topology. These works do not consider an arbitrary splitting of $X$ (Fig.~1(d)); however, such a splitting is very plausible: an agent could have more measurements of some features than of others, and could have measurements that another agent has. In this paper, we fill this gap.

In solving $P_1$, the algorithm we provide is based on the Douglas-Rachford proximal splitting method~\cite{5}. It falls in a growing body of recent work on applying splitting methods, like the alternating directions method of multipliers (ADMM), to obtain distributed algorithms; such work is perhaps most inspired by Boyd et al.~\cite{6}.

\subsection{Contributions}

The contributions of our work are as follows.
\begin{itemize}
\item We present a framework for in-network optimization (Section~II), in which fits $P_1$, and we develop a general algorithm (Algorithm 1).
\item To make $P_1$ amenable to the framework, we provide a new result (Proposition 1) that describes how the constraint in $P_1$ can be separated. 
\item We solve $P_1$ by deriving a specific algorithm (Algorithm 2), and we establish that this algorithm converges.
\end{itemize}

\section{A Framework}

In this section we present a general, variable-centric framework for optimization in a network. Using the framework, we propose a distributed algorithm for solving a general class of problems. The developed framework parallels the strategy described by Parikh and Boyd~\cite{7}, but is specifically adapted to the network setting. It can also be used for the class of problems formalized by Mota et al.~\cite{2}. The tools developed here will be employed later to achieve our goal, to solve $P_1$. 

\subsection{Network Description}

We represent a network by a graph, which consists of nodes and edges. The nodes correspond to agents; the edges, to pairs of agents. As we did for the agents, we label the nodes $1$ to $m$. We let edges represent pairs of agents that can communicate with each other. We denote an edge that joins agents $i$ and $j$ by the ordered pair $(i, j)$, $i < j$. We denote the set of all edges by $\mathcal{E}$, and the set of all nodes forming an edge with $i$ by $\mathcal{N}_i$. The nodes in $\mathcal{N}_i$ are the neighbors of $i$.   

\subsection{Variables}

In general, an optimization problem in a network involves several variables. We assume that the variables are vectors of (possibly different) real finite vector spaces. Some of these variables are specific to each agent, with some neighbors having access to some parts of them.

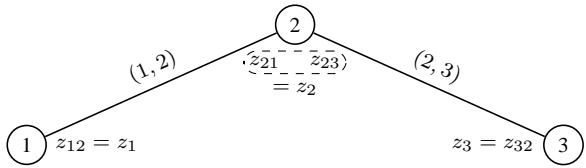
\begin{figure}[t!]
	\centering
\vspace{1.1ex}
	\tikzset{every picture/.style={semithick},font=\footnotesize}
	\begin{tikzpicture}
		\path node at (-3.568,0.4) [circle,draw,inner sep=1mm,label=0:{$z_{12}=z_1$}] (1) {1}
  		      node at ( 0,2) [circle,draw,inner sep=1mm] (2) {2}	
  	       	  node at ( 3.568,0.4) [circle,draw,label=180:{$z_3=z_{32}$},inner sep=1mm] (3) {3}
              node at ( 0,1.499)   [draw,dashed,thin,rounded corners=4.8pt,inner sep=0.7mm,label=-90:{$=z_2$}]{$z_{21}\hspace{3.9mm}z_{23}$} ;	
 
		\draw [-] (1) -- node[above, sloped] {$(1,2)$} (2)
          	  [-] (2) -- node[above, sloped] {$(2,3)$} (3);
	\end{tikzpicture}
	\caption{Variables in a network.}
\end{figure}

Let us consider such variables in a network. As a simple example, we refer to Fig. 2. Each edge involves the variables that can be communicated between agents. Edge $(i, j)$ involves the variable $z_{ij}$. This is the part of the variable of agent $i$ that can be accessed by neighboring agent $j$. And similarly for $z_{ji}$.  

We group all the vectors in the network into a single vector $z$. This vector lies in a real finite vector space $\mathcal{Z}$. We also group the vectors according to edges or nodes, forming two different partitions of $z$. In the example, 
\begin{equation}
z=\underset{\text{\normalsize \raisebox{-1.6mm}{node partition}}}{(z_{12}, (z_{21}, z_{23}), z_{32})}=\underset{\text{\normalsize \raisebox{-1.6mm}{edge partition}}}{((z_{12}, z_{21}), (z_{23}, z_{32})).}
\label{eq:bricks}
\end{equation}
We denote all the variables at each node $i$ by a single variable $z_i$ in a real finite vector space. We refer to this space as $\mathcal{Z}_{i}$. In the example,
\begin{equation}
z_1=z_{12}, \quad z_2=(z_{21},z_{23}), \quad z_3=z_{32},
\label{eq:nodevars}
\end{equation}
and in general, $z_i = (z_{ij})_{j \in \mathcal{N}_i}$.

\subsection{In-Network Optimization}  

The variables that we previously introduced are involved in an optimization problem, and we now give a general description of that problem. 

Let us assign to each node a real-valued cost $g_i(z_{i})$, and to each edge a real-valued cost $g_{ij}(z_{ij},z_{ji})$. We assume that the costs on an edge are the same in both directions, 
\begin{equation*}
g_{ij}(z_{ij},z_{ji})=g_{ji}(z_{ji},z_{ij}),
\end{equation*}
but that they only count as one cost. The idea is to have the agents work together to determine their variables, collectively solving
\begin{equation*}
P_2\colon\quad \min_{z \in \mathcal{Z}}\enskip \sum_{i=1}^m g_i(z_i) + \sum_{(i,j)\in \mathcal{E}}g_{ij}(z_{ij},z_{ji}).
\end{equation*}

\textit{Remark:} Suppose that the agents just needed to minimize the first sum above. They could work independently, each agent $i$ determining $z_i$. They could do this because the objective would be separable, being a sum of terms depending on different parts of $z$, together forming a partition of $z$. Thus, each agent $i$ would just need to minimize $g_i(z_i)$. Similarly, suppose that the agents just needed to minimize the second sum. Pairs of agents could work on separate problems. In $P_2$, however, the variables are all coupled. The parts in the node and edge partitions overlap (see (\ref{eq:bricks}), for example). 

To obtain a distributed approach, we would like the agents to work with each sum in $P_2$ separately. Splitting methods allow us to do this. The most versatile method, requiring only a few assumptions on the problem, is the Douglas-Rachford method~(see [5,~(27.23)], in particular). This method leads to another approach, ADMM~\cite{6}, which has gained popular currency for distributed optimization. Although we can also use ADMM in the present framework, we choose to use the Douglas-Rachford method, because ADMM ends up requiring more exchanges between agents.  

In employing the Douglas-Rachford method to solve $P_2$, the key component is the proximal operator. We denote the proximal operator of a closed proper convex function, say $g\colon \mathcal{Z}\to\mathbf{R}\cup\{+\infty\}$, by $\prox_g$. Applied to a vector in $\mathcal{Z}$, say $\zeta$, the operator gives the unique solution, also in $\mathcal{Z}$, to the following problem:
\begin{equation*} 
\min_{z \in \mathcal{Z}}\enskip g(z) + \frac{1}{2}\lVert z - \zeta \rVert_2^2.
\end{equation*}

Solving our problem using the Douglas-Rachford method involves the proximal operator of each of the sums in the objective of $P_2$. Each operator then entails an optimization problem that can be separated into independent subproblems, with each agent $i$ actually handling $\prox_{g_i}$ and each $(i, j)$ pair handling $\prox_{g_{ij}}$. This decomposition is what allows us to develop an algorithm specifically for the network scenario. We present this algorithm next.

\subsection{In-Network Douglas-Rachford Algorithm}

To solve $P_2$  in a distributed way, we give a description of the part of the Douglas-Rachford method that each agent needs to execute. In that description, we introduce auxiliary variables, $\tilde{z}$ and $\hat{z}$. Note that these variables follow the same indexing scheme and corresponding partition structure as described previously (refer to (\ref{eq:bricks}) and (\ref{eq:nodevars}), for example).

\begin{figure}[h!]
\centering
\fbox{%
\begin{minipage}{7.97cm}
\begin{algorithm}
Choose positive numbers $\lambda$ and $\rho$, with $\rho$ less than $2$. At node $i$ in the network, initialize $z_{i,0}$ to any vector in $\mathcal{Z}_{i}$, and repeat the following: 

After $k$ iterations,
\begin{enumerate}
\item get $z_{ji,k}$ from each neighbor $j$;
\item for each $j$ in $\mathcal{N}_i$, compute $(\tilde{z}_{ij,k+1},\tilde{z}_{ji,k+1})$ from $\prox_{\lambda g_{ij}}(z_{ij,k}, z_{ji,k})$, and assemble $\tilde{z}_{i,k+1}$;
\item $\hat{z}_{i,k+1} = \prox_{\lambda g_i}(2\tilde{z}_{i,k+1}-z_{i,k})$; and
\item $z_{i,k+1}=z_{i,k}+\rho(\hat{z}_{i,k+1}-\tilde{z}_{i,k+1}).$
\end{enumerate}
\end{algorithm} 
\end{minipage}}
\end{figure}

Under certain conditions (see the text following Proposition~2 in the next section, for example), the sequence (of node variables) $\hat{z}_{1}, \hat{z}_{2}, \ldots$ converges to the solution of $P_2$, with each agent having determined its part of that solution.

\section{Solving $P_1$}

The main obstacle in solving $P_1$ is that the constraint couples the data of the agents. Luckily, we can prove the following proposition, which allows us to reformulate the constraint.

\begin{proposition}
There exist in $\mathbf{R}^n$ vectors $\alpha_1$, \dots, $\alpha_m$ that sum to $0$ and for which the set  
\begin{equation*}
\{\beta \in \mathbf{R}^p : \lVert (X_1+\cdots+X_m)\beta - (y_1+\cdots+y_m) \rVert_2 \leq \epsilon \}
\end{equation*}
is equal to the set comprising every $\beta$ such that
\begin{equation*}
 \lVert X_1 \beta - y_1 + \alpha_1 \rVert_2 \leq \frac{\epsilon}{m}, \ldots, \lVert X_m \beta - y_m + \alpha_m \rVert_2 \leq \frac{\epsilon}{m}.
\end{equation*}
\end{proposition}

\subsection{Fitting $P_1$ to the Framework}

By introducing variables $\alpha_1$, \dots, $\alpha_m$ and using Proposition~1, we can fit $P_1$ to the framework described in Section~II.

Let us first specify $z_{ij}$ from the framework to be $(a_{ij}, b_{ij})$, which can be taken as a vector in $\mathbf{R}^{n+p}$. The vector $a_{ij}$ in $\mathbf{R}^n$ is related to $\alpha_i$, and the vector $b_{ij}$ in $\mathbf{R}^p$ is related to $\beta$. To understand these relations, consider the following observation: Since the network is connected, if $b_{ij}=b_{ji}$ along every edge $(i,j)$, and at the same time, at every node $i$, if $b_{ij}$ is equal to some $\beta$ for every $j$ in $\mathcal{N}_i$, then all agents must agree on the same $\beta$. If $\alpha_i = \sum_{j \in \mathcal{N}_i}a_{ij}$, then provided that $a_{ij}=-a_{ji}$ along every edge $(i,j)$, it must be true that $\sum_{i=1}^m\alpha_i=0$.   

We can now describe the edge and node costs. For edge $(i,j)$,
\begin{equation*}
 g_{ij}(z_{ij}, z_{ji}) =
  \begin{cases}
   f(\beta) & \text{if } a_{ij}=-a_{ji} \text{ and } b_{ij}=b_{ji}=\beta;  \\
   +\infty,       & \text{otherwise.}
\end{cases}
\end{equation*}
And for node $i$, 
\begin{equation*}
 g_{i}(z_{i}) =
  \begin{cases}
   0 & \text{if $b_{ij}$ is equal to some $\beta$ for all $j$ in $\mathcal{N}_i$, and}\\
  & \text{$\sum_{j \in \mathcal{N}_i} a_{ij}$ is equal to some $\alpha_i$, with $\beta$ and $\alpha_i$}\\
  & \text{being such that $\lVert X_i\beta-y_i+\alpha_i \rVert_2 \le \epsilon/m$;}\\
     +\infty,       & \text{otherwise.}
  \end{cases}
\end{equation*}
Finally, with these costs, $P_1$ and $P_2$ are equivalent, in the sense that in $P_2$, the $b_{ij}$ part of the solution (the part obtained by agent $i$ for every $j$ in $\mathcal{N}_i$) coincides with the minimizing $\beta$ in $P_1$.

\subsection{In-Network Douglas-Rachford Algorithm for Linear Regression}

Solving $P_1$ is now just a matter of implementing Algorithm~1. Applying the required proximal operators, and simplifying, we obtain the following:
\begin{figure}[h!]
\centering
\fbox{%
\begin{minipage}{7.97cm}
\begin{algorithm}
Choose positive numbers $\lambda$ and $\rho$, with $\rho$ less than $2$. Let $F_\lambda(\cdot)$ denote $2\prox_{(\lambda/2) f}((1/2) \cdot)$. At each node $i$, for all $j$ in $\mathcal{N}_i$, initialize $a_{ij,0}$ to any vector in $\mathbf{R}^n$ and $b_{ij,0}$ to any vector in $\mathbf{R}^p$. Repeat the following at each node $i$: 

After $k$ iterations,
\begin{enumerate}
\item get $a_{ji,k}$ and $b_{ji,k}$ from each neighbor $j$;
\item find the collection $(a_{ij})_{j \in N_i}$ of vectors each in $\mathbf{R}^n$ that add up to some $\alpha_i$ and the vector $\beta$ in $\mathbf{R}^p$ that together minimize
\begin{equation*}
\sum_{j \in \mathcal{N}_i} \Big\{ \lVert a_{ij} +a_{ji,k} \rVert_2^2 + \Vert \beta +b_{ij,k} - F_\lambda (b_{ij,k} + b_{ji,k} ) \rVert_2^2\Big\}
\end{equation*}
subject to the constraint that
\begin{equation*} 
\lVert X_i \beta - y_i + \alpha_i \rVert_2 \leq \frac{\epsilon}{m},
\end{equation*}
and assign the minimizing vectors to $(\hat{a}_{ij, k+1})_{j \in \mathcal{N}_i}$ and $\hat{\beta}_{i,k+1}$; and
\item update:
\begin{equation*}
a_{ij,k+1} = a_{ij,k} -\frac{\rho}{2}(a_{ij,k}-a_{ji,k}) + \rho \hat{a}_{ij,k+1},
\end{equation*}
and
\begin{equation*}
b_{ij,k+1}=b_{ij,k} -\frac{\rho}{2}F_{\lambda}(b_{ij,k}+b_{ji,k}) + \rho\hat{\beta}_{i,k+1}.
\end{equation*}
\end{enumerate}
\end{algorithm} 
\end{minipage}}
\end{figure}

We can say the following about the convergence of Algorithm~2:

\begin{proposition}
At each node $i$, the sequence $\hat{\beta}_{i,1}, \hat{\beta}_{i,2}, \ldots$ converges to the solution of $P_1$ (and thus $P_0$).
\end{proposition}

The proof of Proposition~2 amounts to showing that the two sums in $P_2$ each correspond to closed proper convex functions of $z$, and that those functions satisfy a certain condition (for the particular condition, see [5, Corollary 27.7(a)]).   

\section{A Numerical Experiment}

In this section we provide an example of how our algorithm can be used. 

We consider $P_1$ with $f$ as the $\ell_1$-norm. In this case, $F_\lambda$ in the algorithm is the soft-thresholding operator~[7, (6.9)]. We consider a network of $6$ nodes. The network is generated like this: We randomly pick nodes with replacement. Consecutively picked nodes that are not the same are made neighbors. Keeping track of the picked nodes, we continue this process until all the nodes have been picked. In the network, we consider data from a $20 \times 40$ matrix $X$ and a corresponding vector $y$, both chosen randomly from independent normal entries of mean $0$ and variance $1$. The matrix is split among the agents in the same way as shown in Fig.~1(d). We perform the regression with $\epsilon$ set to $0.01$. For the parameters of the algorithm, we fix $\rho$ to $1.9$ and $\lambda$ to $0.02$. We illustrate in Fig.~3 the convergence of the algorithm. The plot shows the relative error, between the estimate $\hat{\beta}_{i,k}$ of an arbitrary agent $i$ and an estimate, $\hat{\beta}$, computed centrally with $X$ and $y$:
\begin{equation*}
\varepsilon(k) = \frac{\lVert \hat{\beta}_{i,k} -\hat{\beta} \rVert_2}{\lVert \hat{\beta} \rVert_2}.
\end{equation*} 
 The plot depicts the typical error curve~\cite{8} seen when using the Douglas-Rachford method to solve $\ell_1$-minimization problems.

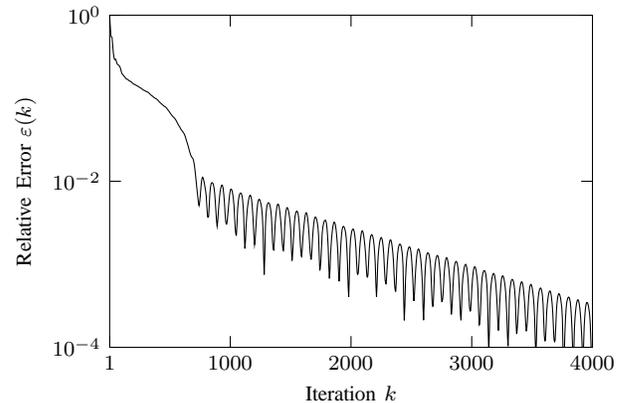
\begin{figure}[h!]  
	\centering 
	\pgfplotsset{compat=1.3,compat=newest,every axis/.append style={font=\footnotesize,thin,tick style={ultra thin, color=black}}} 
	\begin{tikzpicture}
		\begin{semilogyaxis}[%
			xmin=1, xmax=4000,
			ymin=1e-4, ymax=1,
			width=8cm,
			height=6cm,
			ytick={1,1e-2,1e-4},
			xticklabels={$1$,$10^{-2}$,$10^{-4}$},
			xtick={1,1000,2000,3000,4000},
			xticklabels={$1$,$1000$,$2000$,$3000$,$4000$},
			x tick style={color=black, thin},
			y tick style={color=black, thin},
			xlabel={Iteration $k$},
			ylabel={Relative Error $\varepsilon(k)$}
		]
		\addplot[solid,smooth]
			coordinates{(1,0.982471)(11,0.576958)(21,0.531103)(31,0.355108)(41,0.294448)(51,0.294532)(61,0.264049)(71,0.25225)(81,0.245913)(91,0.227483)(101,0.20071)(111,0.191955)(121,0.183807)(131,0.176537)(141,0.169923)(151,0.168835)(161,0.162914)(171,0.160263)(181,0.15718)(191,0.151891)(201,0.147891)(211,0.146143)(221,0.142998)(231,0.140448)(241,0.137976)(251,0.134242)(261,0.132189)(271,0.128593)(281,0.125507)(291,0.123671)(301,0.12133)(311,0.119558)(321,0.117389)(331,0.113915)(341,0.109992)(351,0.107235)(361,0.104469)(371,0.101516)(381,0.0994375)(391,0.0975247)(401,0.0944077)(411,0.0908053)(421,0.0879918)(431,0.085357)(441,0.0827988)(451,0.0805552)(461,0.0788022)(471,0.0767026)(481,0.0731065)(491,0.0701616)(501,0.0666985)(511,0.0636673)(521,0.0609416)(531,0.0589802)(541,0.0569546)(551,0.0540695)(561,0.0506925)(571,0.0474689)(581,0.0451735)(591,0.0423866)(601,0.0405478)(611,0.0386946)(621,0.0365013)(631,0.0331146)(641,0.0295736)(651,0.0266097)(661,0.0233381)(671,0.0209283)(681,0.0195098)(691,0.0187022)(701,0.0161815)(711,0.0122749)(721,0.00883927)(731,0.00648102)(741,0.00511295)(751,0.00656224)(761,0.00984986)(771,0.0112043)(781,0.0101254)(791,0.00843963)(801,0.00654763)(811,0.00386087)(821,0.00390591)(831,0.00718885)(841,0.0092792)(851,0.00961872)(861,0.00898427)(871,0.00708251)(881,0.00435717)(891,0.00288458)(901,0.00424814)(911,0.0065062)(921,0.00825468)(931,0.00901547)(941,0.0081512)(951,0.00569438)(961,0.00340018)(971,0.00307088)(981,0.00417434)(991,0.00619812)(1001,0.00783234)(1011,0.00795978)(1021,0.00687955)(1031,0.0050158)(1041,0.0026885)(1051,0.00228828)(1061,0.00462117)(1071,0.00638769)(1081,0.00714175)(1091,0.00712037)(1101,0.00616819)(1111,0.00386408)(1121,0.0015483)(1131,0.00289286)(1141,0.00459516)(1151,0.00586532)(1161,0.00673062)(1171,0.00656707)(1181,0.00519816)(1191,0.00303514)(1201,0.00172467)(1211,0.00260573)(1221,0.00446991)(1231,0.00576871)(1241,0.00606109)(1251,0.00554284)(1261,0.00456502)(1271,0.00264039)(1281,0.00075207)(1291,0.00284776)(1301,0.00469898)(1311,0.00544028)(1321,0.00546129)(1331,0.00491463)(1341,0.00366884)(1351,0.00179853)(1361,0.00151559)(1371,0.00293306)(1381,0.00434158)(1391,0.00522729)(1401,0.00515907)(1411,0.0041982)(1421,0.00293976)(1431,0.00179028)(1441,0.00149044)(1451,0.00291133)(1461,0.00434086)(1471,0.00482617)(1481,0.00445667)(1491,0.00366191)(1501,0.00251535)(1511,0.00111965)(1521,0.00173495)(1531,0.0032003)(1541,0.00408062)(1551,0.0043729)(1561,0.00406566)(1571,0.00306071)(1581,0.00168674)(1591,0.0012378)(1601,0.00202548)(1611,0.00299191)(1621,0.0039065)(1631,0.00417747)(1641,0.00356521)(1651,0.00246282)(1661,0.00143273)(1671,0.00115729)(1681,0.00201296)(1691,0.00308119)(1701,0.00365916)(1711,0.00363175)(1721,0.00315547)(1731,0.0022405)(1741,0.000877009)(1751,0.00114337)(1761,0.00232249)(1771,0.0030187)(1781,0.00334691)(1791,0.0033442)(1801,0.00273728)(1811,0.00157972)(1821,0.000675797)(1831,0.00135761)(1841,0.00221037)(1851,0.00290798)(1861,0.00322557)(1871,0.00295007)(1881,0.00223247)(1891,0.00136436)(1901,0.000623986)(1911,0.00127248)(1921,0.0022449)(1931,0.00279005)(1941,0.0028372)(1951,0.00257663)(1961,0.0020157)(1971,0.00100641)(1981,0.00040364)(1991,0.00150399)(2001,0.0022489)(2011,0.00260837)(2021,0.0026272)(2031,0.00224635)(2041,0.0015203)(2051,0.000801886)(2061,0.000810498)(2071,0.00146888)(2081,0.00215253)(2091,0.00253007)(2101,0.00238774)(2111,0.00188375)(2121,0.00129408)(2131,0.000735621)(2141,0.000767189)(2151,0.00152553)(2161,0.00209869)(2171,0.00227984)(2181,0.00212803)(2191,0.00169871)(2201,0.000997407)(2211,0.000404828)(2221,0.000989292)(2231,0.00158589)(2241,0.00196461)(2251,0.00210499)(2261,0.00188566)(2271,0.00132398)(2281,0.000689487)(2291,0.000590822)(2301,0.0010154)(2311,0.0015311)(2321,0.00191047)(2331,0.00194184)(2341,0.00162229)(2351,0.0011196)(2361,0.000570149)(2371,0.000485082)(2381,0.00106912)(2391,0.00153252)(2401,0.00173523)(2411,0.0017175)(2421,0.00146607)(2431,0.000931208)(2441,0.000211707)(2451,0.000625789)(2461,0.00115905)(2471,0.00146754)(2481,0.00161583)(2491,0.0015532)(2501,0.00120806)(2511,0.000694079)(2521,0.000362373)(2531,0.000662494)(2541,0.00111106)(2551,0.00144232)(2561,0.00152488)(2571,0.00134738)(2581,0.00101596)(2591,0.000584819)(2601,0.00021254)(2611,0.00066979)(2621,0.00112788)(2631,0.00134823)(2641,0.00136035)(2651,0.0012115)(2661,0.00087304)(2671,0.000371934)(2681,0.000324117)(2691,0.000779401)(2701,0.00109876)(2711,0.00126351)(2721,0.0012452)(2731,0.00101203)(2741,0.000655297)(2751,0.000361911)(2761,0.000409454)(2771,0.000745464)(2781,0.00107059)(2791,0.00120621)(2801,0.00110701)(2811,0.000858797)(2821,0.0005517)(2831,0.00024955)(2841,0.000416596)(2851,0.00078547)(2861,0.00101829)(2871,0.00108326)(2881,0.00100003)(2891,0.000753325)(2901,0.000373222)(2911,0.000212231)(2921,0.000528713)(2931,0.000783904)(2941,0.000951812)(2951,0.000997098)(2961,0.000866784)(2971,0.000591014)(2981,0.000311736)(2991,0.000281068)(3001,0.000516959)(3011,0.000777387)(3021,0.000925388)(3031,0.000905753)(3041,0.000749312)(3051,0.000502614)(3061,0.00019156)(3071,0.00024323)(3081,0.000559843)(3091,0.000760473)(3101,0.000835508)(3111,0.000812329)(3121,0.00066699)(3131,0.000381154)(3141,8.81351e-05)(3151,0.000341148)(3161,0.000571817)(3171,0.000718332)(3181,0.000777447)(3191,0.000715846)(3201,0.000536267)(3211,0.000305167)(3221,0.000159512)(3231,0.00032631)(3241,0.000561641)(3251,0.00070884)(3261,0.000716818)(3271,0.000616609)(3281,0.00045402)(3291,0.00022935)(3301,0.00010358)(3311,0.000363563)(3321,0.000560166)(3331,0.000649263)(3341,0.000648136)(3351,0.000557284)(3361,0.00036916)(3371,0.00014857)(3381,0.000206155)(3391,0.000390119)(3401,0.000533917)(3411,0.000612377)(3421,0.00058639)(3431,0.000458893)(3441,0.000290692)(3451,0.000158689)(3461,0.000207462)(3471,0.000386595)(3481,0.000532003)(3491,0.000573309)(3501,0.000513961)(3511,0.000392338)(3521,0.00022581)(3531,7.94166e-05)(3541,0.000237157)(3551,0.00039937)(3561,0.000492083)(3571,0.00051783)(3581,0.000467901)(3591,0.000329294)(3601,0.000139414)(3611,0.0001251)(3621,0.000270753)(3631,0.000387202)(3641,0.000464239)(3651,0.000471389)(3661,0.000394174)(3671,0.000263068)(3681,0.000128288)(3691,0.000125147)(3701,0.000268694)(3711,0.000391975)(3721,0.000442237)(3731,0.000418902)(3741,0.000342838)(3751,0.000218302)(3761,5.11422e-05)(3771,0.000142435)(3781,0.00028709)(3791,0.000370814)(3801,0.000402286)(3811,0.000382978)(3821,0.00029751)(3831,0.000156573)(3841,6.46369e-05)(3851,0.000175146)(3861,0.000280933)(3871,0.000353829)(3881,0.000373478)(3891,0.000329364)(3901,0.000240664)(3911,0.000135302)(3921,7.00984e-05)(3931,0.000172714)(3941,0.000286428)(3951,0.000343976)(3961,0.000337075)(3971,0.000285361)(3981,0.00019985)(3991,8.21576e-05)(4001,7.43939e-05)(4011,0.00019246)(4021,0.000273956)(4031,0.000313612)(4041,0.000309317)(4051,0.000252393)(4061,0.000151986)(4071,6.40158e-05)(4081,0.000112118)(4091,0.000194675)(4101,0.000263337)(4111,0.000294793)(4121,0.000272214)(4131,0.00020712)(4141,0.000126123)(4151,5.98418e-05)(4161,0.000110482)(4171,0.000200959)(4181,0.000259034)(4191,0.000269003)(4201,0.000239061)(4211,0.000177421)(4221,8.74368e-05)(4231,4.05028e-05)(4241,0.000130187)(4251,0.000199001)(4261,0.000238816)(4271,0.000248212)(4281,0.000216326)(4291,0.00014301)(4301,5.73018e-05)(4311,6.75018e-05)(4321,0.000136057)(4331,0.00019377)(4341,0.000226499)(4351,0.000220998)(4361,0.000180153)(4371,0.000118579)(4381,4.91306e-05)(4391,6.3936e-05)(4401,0.000141495)(4411,0.00019467)(4421,0.000210464)(4431,0.000195748)(4441,0.000156362)(4451,9.04657e-05)(4461,9.2582e-06)(4471,8.0548e-05)(4481,0.000143585)(4491,0.00018106)(4501,0.000194326)(4511,0.00017789)(4521,0.000129376)(4531,6.38662e-05)(4541,3.78651e-05)(4551,8.87552e-05)(4561,0.000140109)(4571,0.000173483)(4581,0.000176763)(4591,0.000150949)(4601,0.000107844)(4611,5.57672e-05)(4621,3.47813e-05)(4631,9.32017e-05)(4641,0.000142911)(4651,0.000164243)(4661,0.000158715)(4671,0.000132195)(4681,8.59551e-05)(4691,2.81885e-05)(4701,4.82544e-05)(4711,9.91699e-05)(4721,0.000134602)(4731,0.000152392)(4741,0.000146441)(4751,0.00011338)(4761,6.38778e-05)(4771,3.02152e-05)(4781,5.81863e-05)(4791,9.87311e-05)(4801,0.000130375)(4811,0.000140679)(4821,0.000126281)(4831,9.45201e-05)(4841,5.32425e-05)(4851,2.08024e-05)(4861,6.08561e-05)(4871,0.000102896)(4881,0.000125173)(4891,0.000126826)(4901,0.000111216)(4911,7.8355e-05)(4921,3.12505e-05)(4931,2.50519e-05)(4941,6.79237e-05)(4951,9.86452e-05)(4961,0.000116441)(4971,0.000117603)(4981,9.76852e-05)(4991,6.12049e-05)};
		\end{semilogyaxis}
	\end{tikzpicture}  
\caption{An example run of Algorithm~2.}
\end{figure}


\begin{thebibliography}{99}
	\bibitem{1} J.~F.~C. Mota, J.~M.~F. Xavier, P.~M.~Q. Aguiar, and M. P\"{u}schel,
``D-ADMM: A communication-efficient distributed algorithm for separable optimization,'' \textit{IEEE Trans. Signal Process.}, vol.~61, no.~10, 2013.
	\bibitem{2} ------, ``Distributed optimization with local domains: Applications in MPC and network flows,'' arXiv:1305.1885, 2013.
	\bibitem{3} ------, ``Distributed compressed sensing algorithms: Completing the puzzle,'' in \textit{Proc. GlobalSIP,} 2013, p.~629.
	\bibitem{4} N. Parikh and S. Boyd, ``Block splitting for distributed optimization,'' \textit{Math. Program. Comput.,} vol.~6, no.~1, pp.~77--102, 2014.
	\bibitem{5} H.~H. Bauschke and P.~L. Combettes, \textit{Convex Analysis and Monotone Operator Theory in Hilbert Spaces.} New York, NY: Springer, 2011.
	\bibitem{6} S. Boyd, N. Parikh, E. Chu, B. Peleato, and J. Eckstein, ``Distributed optimization and statistical learning via the alternating direction method of multipliers,'' \textit{Found. Trends Mach. Learn.,} vol.~3, no.~1, pp. 1--122, 2011.
	\bibitem{7} N. Parikh and S. Boyd, ``Proximal algorithms,'' \textit{Found. Trends in Optim.,} vol.~1, no.~3, pp.~123--231, 2013.
	\bibitem{8} L. Demanet and X. Zhang, ``Eventual linear convergence of the Douglas Rachford iteration for basis pursuit,'' arXiv:1301.0542, 2013. 
\end{thebibliography}
\end{document}